\documentclass[11pt]{amsart}
\usepackage{amssymb, latexsym}
\theoremstyle{plain}
\newtheorem{theorem}{Theorem}

\newtheorem*{thm-cheb}{Theorem (Chebyshev)}

\newtheorem*{2'}{Theorem 2'}
\newtheorem*{3'}{Theorem 3'}

\theoremstyle{remark}

\newtheorem*{Remark 1}{Remark 1}
\newtheorem*{Remark 2}{Remark 2}
\newtheorem*{Remark 3}{Remark 3}
\newtheorem*{Remark 4}{Remark 4}

\numberwithin{equation}{section}

\begin{document}

\title[Kemeny's constant for one-dimensional diffusions]
 {Kemeny's constant for one-dimensional diffusions}

\author{Ross G. Pinsky}


\address{Department of Mathematics\\
Technion---Israel Institute of Technology\\
Haifa, 32000\\ Israel}
\email{ pinsky@math.technion.ac.il}

\urladdr{http://www.math.technion.ac.il/~pinsky/}

\date{}

\begin{abstract}
Let $X(\cdot)$ be a non-degenerate, positive recurrent one-dimensional diffusion process on $\mathbb{R}$ with
 invariant probability
density $\mu(x)$, and
let $\tau_y=\inf\{t\ge0: X(t)=y\}$ denote the first hitting time of $y$.
Let $\mathcal{X}$ be a random variable independent of the diffusion process $X(\cdot)$ and distributed according to the process's invariant probability measure $\mu(x)dx$.
Denote by $\mathcal{E}^\mu$ the expectation with respect to $\mathcal{X}$.
Consider the expression
$$
\mathcal{E}^\mu E_x\tau_\mathcal{X}=\int_{-\infty}^\infty (E_x\tau_y)\mu(y)dy, \ x\in\mathbb{R}.
$$
In words, this expression is the expected hitting time of the diffusion starting from $x$ of a point chosen randomly according to the diffusion's invariant distribution.
We show that this expression is constant in $x$, and that it is finite if and only if
$\pm\infty$ are entrance boundaries for the diffusion.
This result generalizes to diffusion processes the corresponding  result in the setting of  finite Markov chains, where the constant
value is  known as
Kemeny's constant.

\end{abstract}

\maketitle

\section{Introduction and Statement of  Results}
Let $\{X_n\}_{n=0}^\infty$ be an irreducible, discrete time Markov chain on a finite state space $S$, and denote it's invariant
probability measure by $\mu$. For $j\in S$, let  $\hat\tau_j=\inf\{n\ge 1: X_n=j\}$ denote the first passage time to $j$.
Denoting expectations for the process starting from $i\in S$ by $E_i$,
consider the quantity $\sum_{j\in S}\mu_jE_i\hat\tau_j$.
In their book on Markov chains \cite{KS}, Kemeny and Snell showed that the above quantity is independent of the initial
state $i$, and this quantity has become known as \it Kemeny's constant\rm, which we denote by $K$.
Let $\tau_j=\inf\{n\ge 0: X_n=j\}$ denote the first hitting time of $j$. We note that
$\sum_{j\in S}\mu_jE_i\tau_j$ is also independent of $i$, and is equal to $K-1$. This follows from the well-known fact
that  $E_i\hat\tau_i=\frac1{\mu_i}$ \cite{D}.

In \cite{BHLMT}, the authors analysed  the Kemeny constant phenomenon for
 positive recurrent, discrete time and continuous time Markov chains on a  denumerably infinite state space $S$.  They showed
that the quantity $\sum_{j\in S}\mu_jE_i\hat\tau_j$ is either infinite for all $i\in S$, or else is finite and
independent of $i$.
They  conjectured that this quantity
is always infinite in the discrete time setting,  and they proved this in the case of discrete time birth and death chains on $\{0,1,\cdots\}$.
In the case of continuous time birth and   death chains on  $\{0,1,\cdots\}$, they proved that the Kemeny constant is finite if and only if $+\infty$ is an \it entrance boundary\rm\ for the process.

In this paper, we consider the corresponding problem in the context of one-dimensional diffusion processes on $\mathbb{R}$.
Consider a non-degenerate one-dimensional diffusion process $X(\cdot)$ on $\mathbb{R}$ generated by
$$
L=\frac12a(x)\frac{d^2}{dx^2}+b(x)\frac d{dx}.
$$
We  assume that $a$ is continuous and positive, and that $b$ is locally bounded and measurable.
Denote probabilities and expectations for the Markov process $X(\cdot)$ starting from $x\in\mathbb{R}$ by $P_x$ and $E_x$.
For $y\in\mathbb{R}$, let $\tau_y=\inf\{t\ge0: X(t)=y\}$ denote the first hitting time of $y$.
It is well-known \cite{P95} that the following conditions are equivalent:

\begin{equation}\label{posrec}
\begin{aligned}
& i.\ E_x\tau_y<\infty,\ \text{for all}\ x,y\in \mathbb{R};\\
&ii.\ \int_{-\infty}^\infty \frac1{a(x)}\exp\big(2\int_0^x\frac{b(t)}{a(t)}dt\big)dx<\infty;\\
& iii.\ \text{There exists an invariant probability density}\  \mu(x)\ \text{for the process}\ X(\cdot).
\end{aligned}
\end{equation}
If these conditions hold, we say that the process is \it positive recurrent\rm.
In fact then, one has
\begin{equation}\label{invardens}
\mu(x)=\frac {c_0}{a(x)}\exp(2\int_0^x\frac{b(t)}{a(t)})dt,
\end{equation}
 for a normalizing constant $c_0>0$.

From now on we assume that the diffusion is positive recurrent; that is, we assume that
\begin{equation}\label{posrecequ}
\int_{-\infty}^\infty\frac1{a(x)}\exp\big(2\int_0^x\frac{b(t)}{a(t)}dt\big)dx<\infty.
\end{equation}
Let $\mathcal{X}$ be a random variable independent of the diffusion process $X(\cdot)$ and distributed according to the process's invariant probability measure $\mu(x)dx$.
Denote by $\mathcal{E}^\mu$ the expectation with respect to $\mathcal{X}$.
We consider the expression
$$
\mathcal{E}^\mu E_x\tau_\mathcal{X}=\int_{-\infty}^\infty (E_x\tau_y)\mu(y)dy, \ x\in\mathbb{R}.
$$
In words, this expression is the expected hitting time of the diffusion starting from $x$ of a point chosen randomly according to the diffusion's invariant distribution.

There immediately arises the question of whether or not this expression is  finite.
Note the following tradeoff:\it\ On the one hand, the more negative (positive) the drift is in a neighborhood of $+\infty$ ($-\infty$), the faster is the decay of the invariant density $\mu(y)$    at $+\infty$ ($-\infty$). However on the other
hand, the more negative
 (positive) the drift is in a neighborhood of $+\infty$ ($-\infty$), the larger $E_x\tau_y$ will be in a neighborhood of $+\infty$ ($-\infty$).\rm\

  It turns out that the finiteness or infiniteness of the expression depends on whether or not $\pm\infty$ are entrance boundaries for the process.
We recall
 that $+\infty$ is called an \it entrance boundary\rm\
  if \
$\lim_{x\to\infty}P_x(\tau_y<t)>0$, for some $y\in\mathbb{R}$ and some $t>0$.
Similarly, $-\infty$ is called an entrance boundary if
  $\lim_{x\to-\infty}P_x(\tau_y<t)>0$, for some $y\in\mathbb{R}$ and some $t>0$.
  (Actually, equivalently, ``some $y\in\mathbb{R}$ and some $t>0$'' can be replaced by ``all $y\in\mathbb{R}$ and all $t>0$.'')
 Given that the process is positive recurrent, that is, given that \eqref{posrecequ} holds, here is the criterion for an entrance boundary at $+\infty$:
\begin{equation}\label{entranceinfinity}
\int^\infty dx\frac1{a(x)}\exp\big(2\int_0^x\frac{b(s)}{a(s)}ds\big)\int_0^xdy\exp\big(-2\int_0^y\frac{b(s)}{a(s)}ds\big)<\infty.
\end{equation}
See \cite[chapter 8]{P95}, where the term ``explosion inward from infinity'' is used instead of entrance boundary.
The condition \eqref{entranceinfinity} appears as (iv) in Theorem 4.1 in chapter 8. In that theorem, which does not assume positive recurrence, an additional requirement, denoted as (iii),  is also
stated; namely, $\int^{\infty}\exp\big(-2\int_0^x\frac{b(s)}{a(s)}ds\big)dx=\infty$. However, an application of the Cauchy-Schwarz inequality shows that this condition holds automatically if
\eqref{posrecequ} holds.
Similarly,
given that the process is positive recurrent,  here is the criterion for an entrance boundary at $-\infty$:
\begin{equation}\label{entranceminusinfinity}
\int_{-\infty} dx\frac1{a(x)}\exp\big(2\int_0^x\frac{b(s)}{a(s)}ds\big)\int_x^0dy\exp\big(-2\int_0^y\frac{b(s)}{a(s)}ds\big)<\infty.
\end{equation}
We will prove the following theorem. Let $\mu$ denote the probability measure with density $\mu(x)dx$.
\begin{theorem}\label{1}
Assume that the diffusion is positive recurrent; that is, assume that \eqref{posrecequ} holds.
If $\pm\infty$ are both entrance boundaries for the diffusion, that is, if \eqref{entranceinfinity} and \eqref{entranceminusinfinity} both hold, then
$\mathcal{E}^\mu E_x\tau_\mathcal{X}=\int_{-\infty}^\infty (E_x\tau_y)\mu(y)dy$
is finite and independent of $x\in\mathbb{R}$. Two alternative expressions for the value of this constant are
\begin{equation}\label{constant}
2\int_{-\infty}^\infty dy\thinspace \mu(y)\int_y^\infty dz\frac{\mu([z,\infty))}{\mu(z)a(z)}\ \ \text{and}\ \
2\int_{-\infty}^\infty dy\thinspace \mu(y)\int_{-\infty}^y dz\frac{\mu((-\infty,z])}{\mu(z)a(z)}.
\end{equation}
If at least one  of $\pm\infty$ is not an entrance boundary, that is if at least one of \eqref{entranceinfinity} and \eqref{entranceminusinfinity} does not hold, then
\begin{equation}
\mathcal{E}^\mu E_x\tau_\mathcal{X}=\int_{-\infty}^\infty (E_x\tau_y)\mu(y)dy=\infty,\ \text{for all}\ x\in\mathbb{R}.
\end{equation}
\end{theorem}
\bf\noindent Remark 1.\rm\ Given a continuously differentiable, strictly positive probability density $\mu$ and given a continuously differentiable diffusion matrix $a$,
if one chooses the drift $b(x)=\frac12\big(a(x)\frac{\mu'(x)}{\mu(x)}+a'(x)\big)$,
then the diffusion process with generator  $L$ will have invariant probability density $\mu$.
Thus, given such a density $\mu$, the diffusion processes for which $\mu$ is the invariant density can be
indexed by their diffusion matrices $a$. From \eqref{constant} we see that given the invariant density $\mu$, the expression
$\mathcal{E}^\mu E_x\tau_\mathcal{X}$ is monotone decreasing  as a function of the diffusion matrix $a$.
Furthermore, we see that for sufficiently large $a$ it will be finite and for sufficiently small $a$ it will be infinite.
In particular then, given $\mu$ we can find a diffusion with invariant  density $\mu$ for which $\pm\infty$ are entrance boundaries and we can find such a diffusion for which
$\pm\infty$  
are not entrance boundaries. 
\medskip

\noindent \bf Remark 2.\rm\ Let $\mu$ be a continuously differentiable, strictly positive probability density as in  Remark 1. 
Since the two expressions in \eqref{constant} must be both  finite or both infinite, it is easy to see that in the case of constant diffusion coefficient, $a\equiv \text{const.}$,
the expression 
$\mathcal{E}^\mu E_x\tau_\mathcal{X}$  is finite if and only if 
\begin{equation}\label{condition}
\int_{-\infty}^\infty\frac{\mu([y,\infty))}{\mu(y)}dy<\infty.
\end{equation}
In particular, if $\mu(x)\sim\text{const.}e^{-k|x|^l}$, for $k,l>0$, then \eqref{condition} holds if and only if $l>2$,

\section{Proof of Theorem \ref{1}}
We first proof that $\int_{-\infty}^\infty (E_x\tau_y)\mu(y)dy<\infty$
if and only if \eqref{entranceinfinity} and \eqref{entranceminusinfinity} hold.
We have the following explicit expression for the expected hitting time:
\begin{equation}\label{exphitting}
E_xT_y=\begin{cases}
2\int_y^xdz\exp(-2\int_0^z\frac{b(t)}{a(t)}dt)\int_z^\infty dw\frac1{a(w)}\exp(2\int_0^w\frac{b(t)}{a(t)}dt),\ -\infty< y<x;\\
 2\int_x^ydz\exp(-2\int_0^z\frac{b(t)}{a(t)}dt)\int_{-\infty}^zdw\frac1{a(w)}\exp(2\int_0^w\frac{b(t)}{a(t)}dt),\ x<y<\infty.
\end{cases}
\end{equation}
For a derivation, see for example the proof of Proposition 2 in \cite{P} (where $a(x)$ is a constant and denoted by $D$).
Using this with \eqref{invardens}--\eqref{entranceminusinfinity}, it is easy to see that
\eqref{entranceinfinity} and \eqref{entranceminusinfinity} constitute necessary and sufficient conditions for
the finiteness of $\int_{-\infty}^\infty (E_x\tau_y)\mu(y)dy$. Indeed, from \eqref{exphitting} and \eqref{invardens}, we have
\begin{equation}\label{keyforfinite}
\begin{aligned}
&\int_x^\infty(E_x\tau_y)\mu(y)dy=\\
&2\int_x^\infty  dy\frac {c_0}{a(y)}\exp\big(2\int_0^y\frac{b(t)}{a(t)}dt\big) \int_x^ydz\exp(-2\int_0^z\frac{b(t)}{a(t)}dt)\int_{-\infty}^zdw\frac1{a(w)}\exp(2\int_0^w\frac{b(t)}{a(t)}dt).
\end{aligned}
\end{equation}
By \eqref{invardens}, the right hand side of \eqref{keyforfinite} is finite if and only if
$$
\int_x^\infty  dy\frac {c_0}{a(y)}\exp\big(2\int_0^y\frac{b(t)}{a(t)}dt\big) \int_x^ydz\exp(-2\int_0^z\frac{b(t)}{a(t)}dt)<\infty,
$$
and  this latter expression is finite if and only if \eqref{entranceinfinity} holds.
Thus, $\int_x^\infty(E_x\tau_y)\mu(y)dy<\infty$ if and only if \eqref{entranceinfinity} holds.
A similar analysis shows that $\int_{-\infty}^x(E_x\tau_y)\mu(y)dy<\infty$ if and only if
\eqref{entranceminusinfinity} holds.

We now show that if \eqref{entranceinfinity} and \eqref{entranceminusinfinity} hold, then $\int_{-\infty}^\infty(E_x\tau_y)\mu(y)dy$ is independent of $x\in\mathbb{R}$.
For $y\in \mathbb{R}$, define $u_{y,+}(x)=E_x\tau_y$, for $x\ge y$, and  $u_{y,-}(x)=E_x\tau_y$, for $x\le y$.
(Of course, $u_{y,+}(x)$ and $u_{y,-}(x)$ respectively are equal to the first and second lines on the right hand side of \eqref{exphitting}.)
As is well-know, it follows from an application of Ito's formula that
\begin{equation}\label{pderep}
\begin{aligned}
&Lu_{y,+}=-1\ \text{in}\ (y,\infty);\ u_{y,+}(y)=0;\\
&Lu_{y,-}=-1\ \text{in}\ (-\infty, y);\ u_{y,-}(y)=0.
\end{aligned}
\end{equation}
(Indeed, it is from this that the formulas in \eqref{exphitting} were derived.)
Define $F(x)=\int_{-\infty}^\infty(E_x\tau_y)\mu(y)dy$. Then we have
$$
F(x)=\int_{-\infty}^x u_{y,+}(x)\mu(y)dy+\int_x^\infty u_{y,-}(x)\mu(y)dy.
$$
In light of the fact that  $u_{y,+}(x)$ and $u_{y,-}(x)$ are given by  \eqref{exphitting}, as well as the fact that \eqref{posrecequ}-- \eqref{entranceminusinfinity} hold,
we can differentiate freely under the integral.
Using the boundary condition in \eqref{pderep}, we have
\begin{equation}\label{firstderiv}
F'(x)=\int_{-\infty}^x u'_{y,+}(x)\mu(y)dy+\int_x^\infty u'_{y,-}(x)\mu(y)dy,
\end{equation}
and differentiating again gives
\begin{equation}\label{secondderiv}
F''(x)=\int_{-\infty}^x u''_{y,+}(x)\mu(y)dy+\int_x^\infty u''_{y,-}(x)\mu(y)dy+\mu(x)\big(u_{x,+}'(x)-u_{x,-}'(x)\big).
\end{equation}
From \eqref{firstderiv} and \eqref{secondderiv} we obtain
\begin{equation}\label{LF}
LF(x)=\int_{-\infty}^x \mu(y)Lu_{y,+}(x)dy+\int_x^\infty\mu(y)Lu_{y,-}(x)dy+\frac12a(x)\mu(x)\big(u_{x,+}'(x)-u_{x,-}'(x)\big).
\end{equation}

From \eqref{pderep} we have
\begin{equation}\label{firstpart}
\int_{-\infty}^x \mu(y)Lu_{y,+}(x)dy+\int_x^\infty\mu(y)Lu_{y,-}(x)dy=-\int_{-\infty}^\infty\mu(y)dy=-1.
\end{equation}
Using the formulas for  $u_{y,+}(x)$ and $u_{y,-}(x)$ as given by the two lines on the right hand side of \eqref{exphitting}, and recalling \eqref{invardens},
we have
$$
\begin{aligned}
&u_{x,+}'(x)-u_{x,-}'(x)=
2\exp(-2\int_0^x\frac{b(t)}{a(t)}dt)\int_x^\infty dw\frac1{a(w)}\exp(2\int_0^w\frac{b(t)}{a(t)}dt)-\\
&2\exp(-2\int_0^x\frac{b(t)}{a(t)}dt)\int_{-\infty}^xdw\frac1{a(w)}\exp(2\int_0^w\frac{b(t)}{a(t)}dt)=\\
&2\exp(-2\int_0^x\frac{b(t)}{a(t)}dt)\int_{-\infty}^\infty dw\frac1{a(w)}\exp(2\int_0^w\frac{b(t)}{a(t)}dt)=\frac2{c_0}
\exp(-2\int_0^x\frac{b(t)}{a(t)}dt).
\end{aligned}
$$
Thus,
\begin{equation}\label{secondpart}
\frac12a(x)\mu(x)\big(u_{x,+}'(x)-u_{x,-}'(x)\big)=\frac{c_0}2\exp(2\int_0^x\frac{b(t)}{a(t)}dt)\times
\frac2{c_0}
\exp(-2\int_0^x\frac{b(t)}{a(t)}dt)=1.
\end{equation}
From \eqref{LF}--\eqref{secondpart}, we conclude that $LF=0$; that is, $F$ is $L$-harmonic.

Since $L$ is a recurrent diffusion generator, it has no nonconstant positive harmonic functions \cite[p.457]{P95}. Consequently,
we conclude that $F(x)=\int_{-\infty}^\infty(E_x\tau_y)\mu(y)dy$ is constant in $x$.

It remains to prove \eqref{constant}.
From \eqref{keyforfinite} and the corresponding formula for $\int_{-\infty}^x(E_x\tau_y)\mu(y)dy$,
we have
\begin{equation}\label{final}
\begin{aligned}
&F(x)=\int_{-\infty}^\infty(E_x\tau_y)\mu(y)dy=\\
&2\int_x^\infty  dy\thinspace\mu(y) \int_x^ydz\exp(-2\int_0^z\frac{b(t)}{a(t)}dt)\int_{-\infty}^zdw\frac1{a(w)}\exp(2\int_0^w\frac{b(t)}{a(t)}dt)+\\
&2\int_{-\infty}^x  dy\thinspace \mu(y) \int_y^xdz\exp(-2\int_0^z\frac{b(t)}{a(t)}dt)\int_z^\infty dw\frac1{a(w)}\exp(2\int_0^w\frac{b(t)}{a(t)}dt).
\end{aligned}
\end{equation}
Letting $x\to\infty$ in \eqref{final} and using \eqref{invardens} to write everything in terms
of $a$ and $\mu$ gives the first alternative in \eqref{constant}. Similarly, letting
$x\to-\infty$ gives the second alternative in \eqref{constant}.\hfill$\square$

\end{document}